\documentclass[11pt]{article}
\usepackage{amsfonts,amsmath,amssymb}

\newtheorem{lemma}{Lemma}

\newcommand{\Mc}{\mathcal{M}}
\newtheorem{theorem}{Theorem}
\newtheorem{corollary}{Corollary}
\newtheorem{remark}{Remark}

\newcommand{\R}{\mathbb{R}}
\newcommand{\Tb}{\mathbb{T}}
\newcommand{\N}{\mathbb{N}}

\newcommand{\Z}{\mathbb{Z}}

\author{
 
 Yuri Bakhtin\footnotemark[1],
 Christine Heitsch\footnotemark[2] 
}

\title{Large Deviations for Random Trees}

\begin{document}
\maketitle
\footnotetext[1]{School of Mathematics, Georgia Tech, Atlanta GA, 30332-0160; email:~bakhtin@math.gatech.edu,}
\footnotetext[2]{School of Mathematics, Georgia Tech, Atlanta GA, 30332-0160; email:~heitsch@math.gatech.edu}
	
\abstract{We consider large random trees under Gibbs distributions and prove a Large Deviation Principle~(LDP) for the distribution of degrees of vertices of the tree. The LDP rate function is given explicitly. An immediate consequence
is a Law of Large Numbers for the distribution of vertex degrees in a large random tree.
Our motivation for this study comes from the analysis of RNA secondary structures.

} 

\section{Introduction}
In this note, we prove a Large Deviation Principle (LDP) for two models of
equilibrium statistical mechanics.
In both cases, we consider a set of trees on $N$ vertices and
we define the Gibbs distribution associated to a certain energy
function on that set.
The main goal of our work is to study some typical features of large
random trees ($N\to\infty$) under these distributions.

Here, we provide rigorous proofs for the LDP results announced
in~\cite{yb-ceh-bmb}.  As discussed there, our results are motivated
by, and have applications to, the branching of RNA secondary structures.
The trees we consider are a useful abstraction of these biological
structures (see~\cite{plane, insights} for references on this connection)
as well as relatively straightforward to analyze mathematically.
In this simplified model of RNA folding, we can address the interplay
between entropy and energy in determining a ``typical'' branching
configuration.
We find that, due to the entropy factor, the typical
configurations in our model differ from the arrangements which have
minimal energy in interesting ways. 

Our mathematical results support and extend recent developments in
RNA secondary structure prediction
(reviewed in~\cite{mathews-06, mathews-turner-06})
which broaden the focus beyond simply finding a structure with
minimal free energy.
In particular, we prove a Law of Large Numbers for the degree frequencies
in our large random trees, and find that the most common trees are not
the minimizers of the associated energies.
This highlights the limitations of prediction methods focused solely
on energy minimization and the significance of entropy considerations
in computational structural biology.

\section{Models and results} 
In this section we describe our models
and state the results.  The proofs are given in the next section.
\subsection{Labeled trees}\label{subsec:labeled}
In our first model we fix a natural number $D\ge 2$ and for each $N\in\N$ consider the set $\Tb_N(D)$ of labeled trees on $N\in\N$ vertices such that the
degree of each vertex does not exceed $D$. To define Gibbs distributions on~$\Tb_N(D)$ we need a 
function $c:\{1,\ldots,D\}\to\R$ which plays the role of the energy associated with the degree of a vertex.

To each of the trees $T$ in $\Tb_N(D)$ we associate the energy
\begin{equation}
H(T)=\sum_{j=1}^N c(d_j(T))=\sum_{k=1}^D c(k)\chi_k(T),
\end{equation}
where $d_j(T)$ denotes the degree of the $j$-th vertex,
and $\chi_k(T)$ is the number of vertices of degree $k$ in $T$. 
Now the Gibbs probability measure on $\Tb_N(D)$ associated with $H$ is given by
\begin{equation*}
P_N\{T\}=\frac{e^{-\beta H(T)}}{Z_N},\quad T\in\Tb_N(D),
\end{equation*}
where $\beta>0$ is the  inverse temperature parameter and 
\begin{equation}
\label{eq:Z_N}
Z_N=\sum_{T\in\Tb_N} e^{-\beta H(T)}
\end{equation}
is the partition function.

Our first result is an LDP for the
degree distribution of random labeled trees under measures $P_N$ introduced above.

Let us recall that a sequence of probability measures $(\mu_N)_{N\in\N}$ on a compact metric space $(E,\rho)$ satisfies an LDP with a lower-semicontinuous 
nonnegative rate function $I:E\to\R$ if  
\begin{equation*}
\limsup_{N\to\infty}\frac{1}{N} \ln\mu_N(C)\le -I(C),\quad \mbox{\rm for any closed set\ }C\subset E,
\end{equation*}
and
\begin{equation*}
\liminf_{N\to\infty}\frac{1}{N} \ln\mu_N(O)\ge -I(O),\quad \mbox{\rm for any open set\ }O\subset E,
\end{equation*}
where for $U\subset E$,
\begin{equation*}
I(U)=\inf_{p\in U} I(p).
\end{equation*}
See \cite[Section II.3]{Ellis:MR2189669} or \cite[Section 1.2]{Dembo-Zeitouni:MR1619036} for further details.

Informally, an LDP means that if we consider random variables $X_N$ with distribution $\mu_{N}$, then for all $p$ and large $N$
we have $$\mu_N\{X_N\approx p\}\approx e^{- NI(p)}.$$ In particular, if the minimal value $0$ is attained by $I$ at a unique point~$p^*$
then for any neighborhood $O$ of $p^*$, $\mu_N(O^c)$ decays exponentially in~$N$. This can be restated as a Law of Large Numbers with exponential convergence
in probability to the limit point $p^*$.    

We can view $(\chi_1,\ldots,\chi_D)$ as
a random vector defined on the probability space $\Tb_N(D)$ equipped with the Gibbs measure $P_N$. We would like to study the frequencies of vertex degrees, so for each $N$ we
introduce a probability measure $\nu_N$ on $[0,1]^D$ defined as the distribution
of the random vector $\frac1N (\chi_1,\ldots,\chi_D)$ under $P_N$.
It is natural to formulate an LDP for $\nu_N$ on the set
\begin{equation*}
\Mc=\left\{p\in[0,1]^D:\ \sum_{k=1}^Dp_k=1,\ \sum_{k=1}^D k p_k=2\right\}
\end{equation*}
equipped with Euclidean distance. (Notice that $\Mc$ is nonempty if $D\ge 2$.)
Though the random vector $\frac1N (\chi_1,\ldots,\chi_D)$ does not belong to $\Mc$,
it is asymptotically close to $\Mc$:
\begin{equation*}
\sum_{k=1}^D\frac{\chi_k}{N}=1,\quad \sum_{k=1}^Dk\frac{\chi_k}{N}=2-\frac{2}{N}.
\end{equation*}
So instead of formulating an LDP for the sequence of random vectors $\frac1N (\chi_1,\ldots,\chi_D)$,
we shall formulate and prove an LDP for a sequence of random vectors that is close to it and belongs to $\Mc$.

To define the rate function, we introduce $J:\Mc\to\R$ via
\begin{equation*}
J(p)=-h(p)+\beta E(p)+G(p),
\end{equation*}
where
\begin{equation*}
h(p)=-\sum_{k=1}^D p_k\ln p_k
\end{equation*}
is the entropy of the probabilty vector $p=(p_1,\ldots,p_D)$,
\begin{equation*}
E(p)=\sum_{k=1}^D p_kc(k)
\end{equation*}
is the energy associated with $p$, and $G(p)$ is defined by
\begin{equation}
\label{eq:G}
G(p)=\sum_{k=1}^{D}p_k\ln ((k-1)!).
\end{equation}
In Section~\ref{sec:proof}, we shall see that
the function $G$ appears naturally in the analysis of random trees. 

The function $J$ is strictly convex down and continuous on $\Mc$. Therefore,
it attains its minimal value at a uniquely defined point $p^*\in\Mc$.
Consider now
\begin{equation}
\label{eq:rate}
I(p)=J(p)-J(p^*).
\end{equation} 
It is easy to see that $I$ is bounded, convex and continuous on $\Mc$.

For a measure $Q$ on $[0,1]^D\times\Mc$ we define  $Q^{(1)}$ and $Q^{(2)}$
as the marginal distributions of $Q$ on $[0,1]^D$ and $\Mc$ respectively. 

\begin{theorem}\label{th:mainLDP} There is a sequence of probability measures $(Q_N)_{N\in\N}$ defined on
$[0,1]^D\times\Mc$ with the following properties.
\begin{enumerate}
\item For each $N$, we have $Q_N^{(1)}=\nu_N$.
\item For each $N$, 
\begin{equation*}
Q_N\left\{(x,y)\in[0,1]^D\times\Mc:\sum_{k=1}^D|x_k-y_k|>\frac{2}{N}\right\}=0.
\end{equation*}
\item The sequence $(Q_N^{(2)})_{N\in\N}$ satisfies an LDP on $\Mc$ with the rate function $I$ defined
in \eqref{eq:rate}.
\end{enumerate}
\end{theorem}

\begin{remark}\rm This theorem says that although the random vector $\chi/N$ does not belong to $\Mc$,
one can find another random vector that is, on the one hand, very close to $\chi/N$ and on the other hand
belongs to $\Mc$ and satisfies the LDP.
\end{remark}

Theorem~\ref{th:mainLDP} immediately implies the following Law of Large Numbers: 
\begin{corollary} As $N\to\infty$,
\begin{equation*} 
\left(\frac{\chi_1}{N},\ldots, \frac{\chi_D}{N}\right)\to p^*
\end{equation*}
in probability.
\end{corollary}
\begin{remark}\rm
The statements above show that with high probability the degree
frequencies are close to $p^*$. Note that in most cases the minimum of 
the energy $E$ on $\Mc$ is not attained at $p^*$.
\end{remark}

\subsection{Plane trees}\label{subsec:plane}
We now consider a similar model for plane trees (sometimes also called ordered trees). These are rooted trees
such that subtrees at any vertex are linearly ordered, see e.g.~\cite{Stanley:MR1676282}. We redefine the notation
introduced in the previous section. We fix a number $D\in\N$ and for each $N\in\N$ let $\Tb_N(D)$ denote the set of ordered trees on $N\in\N$ vertices such that the
branching (i.e.\ the number of children) at each vertex does not exceed $D$.  The energy
of each vertex depends only on its branching and is given by
a function $c:\{0,1,\ldots,D\}\to\R$. With each tree $T\in\Tb_N(D)$ we associate the energy
\begin{equation}
H(T)=\sum_{k=0}^D c(k) \chi_k(T),
\end{equation}
where $\chi_k(T)$ is now the number of vertices with $k$ children in $T$. 
The Gibbs probability measure on $\Tb_N(D)$ associated with $H$ is given by
\begin{equation*}
P_N\{T\}=\frac{e^{-\beta H(T)}}{Z_N},\quad T\in\Tb_N(D),
\end{equation*}
where $\beta>0$ is the inverse temperature and 
$Z_N$ is a normalizing constant.

For each $N$, we introduce 
a probability measure $\nu_N$ on $[0,1]^{D+1}$ defined as the distribution of the random vector $\frac1N (\chi_0, \chi_1,\ldots,\chi_D)$ under $P_N$.

We redefine $\Mc$ to be
\begin{equation*}
\Mc=\left\{p\in[0,1]^{D+1}:\ \sum_{k=0}^Dp_k=1,\ \sum_{k=0}^D k p_k=1\right\}.
\end{equation*}
To formulate an LDP for this model we define 
$J:\Mc\to\R$ via
\begin{equation*}
J(p)=-h(p)+\beta E(p),
\end{equation*}
where
\begin{equation*}
h(p)=-\sum_{k=0}^D p_k\ln p_k
\end{equation*}
is the entropy of the probabilty vector $p=(p_0,p_1\ldots,p_D)$,
and 
\begin{equation*}
E(p)=\sum_{k=0}^D p_kc(k)
\end{equation*}
is the energy associated with $p\in\Mc$.

As in the first model, the function $J$ attains
its minimum on $\Mc$ at a unique point that we denote by $p^*$.
Let
\begin{equation}
\label{eq:rate_ordered}
I(p)=J(p)-J(p^*).
\end{equation}
This function will play the role of the rate function. Notice that in the case of plane trees it does not involve
the function $G(p)$ that appeared in the construction of the rate function for the case of labeled trees.

For a measure $Q$ on $[0,1]^{D+1}\times\Mc$ we define  $Q^{(1)}$ and $Q^{(2)}$
as the marginal distributions of $Q$ on $[0,1]^{D+1}$ and $\Mc$ respectively. 

\begin{theorem}\label{th:mainLDP-ordered} 
 There is a sequence of probability measures $(Q_N)_{N\in\N}$ defined on
$[0,1]^{D+1}\times\Mc$ with the following properties.
\begin{enumerate}
\item For each $N$, we have $Q_N^{(1)}=\nu_N$.
\item For each $N$, 
\begin{equation*}
Q_N\left\{(x,y)\in[0,1]^{D+1}\times\Mc:\sum_{k=0}^D|x_k-y_k|>\frac{1}{N}\right\}=0.
\end{equation*}
\item The sequence $(Q_N^{(2)})_{N\in\N}$ satisfies an LDP on $\Mc$ with the rate function $I$ defined
in \eqref{eq:rate_ordered}.
\end{enumerate}
\end{theorem}
An immediate consequence is the following Law of Large Numbers:
\begin{corollary}\label{cr:LLN2} As $N\to\infty$,
\begin{equation*} 
\left(\frac{\chi_0}{N},\frac{\chi_1}{N},\ldots, \frac{\chi_D}{N}\right)\to p^*
\end{equation*}
in probability.
\end{corollary}

\section{Proofs}\label{sec:proof}

We start with the proof of Theorem~\ref{th:mainLDP}, adopting the notation and setting for labeled trees from Section~\ref{subsec:labeled}.

The crucial fact for our analysis is the following formula for
the number of trees on $N$ vertices with degrees $d_1,\ldots,d_N$:
\begin{equation*}
\binom{N-2}{d_1-1,\ d_2-1,\ \ldots,\ d_N-1 }
\end{equation*}
if $d_1+\ldots+d_N=2N-2$, and $0$ otherwise,
see \cite[Formula (2.1)]{Moon:MR0274333}. Therefore, the total number of $N$-trees $T$ with $\chi(T)=(n_1,\ldots,n_D)$
is given by
\begin{multline*} 
\binom{N-2}{\underbrace{0,\ldots,0}_{n_1}, \underbrace{1,\ldots,1}_{n_2}, \ldots, \underbrace{D-1,\ldots,D-1}_{n_D}}\binom{N}{n_1,\ldots,n_D}
\\ =\frac{(N-2)!}{(2!)^{n_3}\ldots ((D-1)!)^{n_D}}C(N,n),
\end{multline*}
where $C(N,n)=\binom{N}{n_1,\ldots,n_D}$.
All these trees $T$ have the same energy $H(T)$, so that
\begin{equation}
P_N\left\{\frac{\chi(T)}{N}=\frac{n}{N}\right\}= \frac{e^{-N F\left(\frac{n}{N}\right)}C(N,n)}{
Z_N},
\label{eq:projection_of_P_N}
\end{equation} 
where $Z_N$ is defined in~\eqref{eq:Z_N}, and we notice that
\begin{equation*}
Z_N= \sum_{\substack{n_1+\ldots+n_D=N\\ n_1+\ldots+Dn_D=2N-2}} e^{-N F\left(\frac{n}{N}\right)}C(N,n),
\end{equation*} 
and
\begin{equation*}
F(p)=\beta E(p)+G(p)=\beta \sum_{k=1}^Dc(k)p_k+\sum_{k=1}^D \ln((k-1)!)p_k,\quad p\in[0,1]^D,
\end{equation*}
with $G(p)$ defined in~\eqref{eq:G}.

Our plan is to use the LDP for multinomial distribution that manifests itself
in coefficients $C(N,n)$ in the r.h.s. of \eqref{eq:projection_of_P_N}, and then apply a version of Varadhan's lemma for
Gibbs transformation via the exponential factor $e^{-N F\left(\frac{n}{N}\right)}$.
 
We start with
the family of distributions $\mu_N$ on $\Mc$ defined by 
\begin{equation*}
\mu_N\left\{\left(\frac{n_1}{N},\ldots,\frac{n_D}{N}\right)\right\}= \begin{cases}\frac{C(N,n)}{ Z_N'},&\mbox{if\ } \left(\frac{n_1}{N},\ldots,\frac{n_D}{N}\right)\in\Mc \\0,&\mbox{\rm otherwise} \end{cases},
\end{equation*}
where
\begin{equation*}
Z'_N=\sum_{n/N\in\Mc} C(N,n).
\end{equation*}

\begin{lemma} The sequence of measures $(\mu_N)_{N\in\N}$ satisfies an LDP on $\Mc$ with rate function $I_1$ defined by $$I_1(p)=h^*-h(p),$$
where $$h^*=\sup_{p\in\Mc} h(p).$$
\end{lemma}
Proof. The proof of this lemma literally repeats that of Sanov's theorem (an LDP for the multinomial distribution, see \cite[Theorem 2.1.10]{Dembo-Zeitouni:MR1619036}). It is based on the formula:
\begin{equation*}
\frac{1}{N}\ln C(N,n)= -\sum_{k=1}^D\frac{n_k}{N}\ln \frac{n_k}{N}+ O\left(\frac{\ln N}{N}\right), \mbox{\rm\ as\ } N\to\infty,
\end{equation*}
which holds true uniformly in $n$, see e.g.\cite[Lemma I.4.4]{Ellis:MR2189669}.

Let us now introduce the Gibbsian weight
\begin{equation*}
q_N\left(\frac{n}{N}\right)=e^{-N F\left(\frac{n}{N}\right)},
\end{equation*}
and a new family of measures $\lambda_N$ on $\Mc$:
\begin{equation*}
\lambda_N\left\{\frac{n}{N}\right\}=\frac{q_N\left(\frac{n}{N}\right)\mu_N
\left\{\frac{n}{N}\right\}}{Z''_N},\quad\text{for}\ \frac{n}{N}\in\Mc, 
\end{equation*}
where
\begin{equation*}
Z''_N=\sum_{\frac{n}{N}\in\Mc}q_N\left(\frac{n}{N}\right)\mu_N\left\{\frac{n}{N}\right\}
=\int_{\Mc} e^{-NF(p)}\mu_N(dp).
\end{equation*}
In other words,
\begin{equation*}
\lambda_N(dp)=\frac{e^{-NF(p)}\mu_N(dp)}{\int_{\Mc} e^{-NF(p)}\mu_N(dp)}.
\end{equation*}
Let us also denote $J_1(p)=F(p)+I_1(p)$ and $J_{1,*}=\inf_{p\in\Mc}J_1(p)$.
\begin{lemma}\label{lm:LDP-aux} The sequence of measures $(\lambda_N)_{N\in\N}$ satisfies an LDP on $\Mc$ with rate function $I_2$ given by
$I_2(p)=J_1(p)-J_{1,*}$.
\end{lemma}
Proof. This lemma follows directly from a variant of Varadhan's lemma for Gibbs transformations (Theorem II.7.2 in [Ellis]).
\begin{remark} \rm
Notice that $I_2(p)=I(p)$ for all $p\in\Mc$. So we have proven the desired LDP on $\Mc$ for $(\lambda_N)_{N\in\N}$, and in order to prove 
Theorem~\ref{th:mainLDP} we shall have to compare $\lambda_N$ to $\nu_N$.
\end{remark}
\medskip

\noindent
Proof of Theorem~\ref{th:mainLDP}. We consider the distribution $P_N$ on $\Tb_N(D)$, so that $\frac{\chi}{N}$ is distributed according to $\nu_N$.
For each $x$ that belongs to the support of $\nu_N$ we introduce the set
\begin{multline*}
R(x)=\Biggl\{y\in\Mc:\ y_k=\frac{m_k}{N}, m_k\in\Z, k=1,\ldots,D,\\ \mbox{\rm and}\ \sum_{k=1}^D|x_k-y_k|=\frac 2N\Biggr\}. 
\end{multline*}
It is easy to see that $1\le|R(x)|\le D^2$
for all $x$, where $|R|$ denotes the number of elements in $R$.

Let us now define the measure $Q_N$. We start with random variables $\chi/N$, and define
a random vector $Y$ so that, given
$\chi/N$, the conditional distribution of $Y$ is uniform on $R(\chi/N)$. Now $Q_N$
denotes the joint distribution of $\chi/N$ and $Y$. Clearly, the first two desired 
properties of $Q$ hold true by the definition of~$Q_N$. The third one follows from
Lemma~\ref{lm:LDP-aux} and the following statement claiming that 
measures $Q_N^{(2)}$ and $\lambda_N$ differ by a subexponential factor,  thus
obeying an LDP with the same rate function: 
\begin{lemma} There is a constant $C>0$ such that for all $N$
and all sets $U\subset\Mc$,
\begin{equation*}
\frac{1}{CN^4}\le \frac{Q_N^{(2)}(U)}{\lambda_N(U)}\le CN^4.
\end{equation*}
\end{lemma}
This lemma is a straightforward consequence of the following fact: there is a constant $K$ such that
if $|n_1-n'_1|+\ldots+|n_D-n'_D|=2$ then
$$
\frac{1}{K N^2}\le \frac{e^{-N F\left(\frac{n}{N}\right)}C(N,n)}{e^{-N F\left(\frac{n'}{N}\right)}C(N,n')}\le K N^2.
$$

\bigskip

The proof of Theorem~\ref{th:mainLDP-ordered} is essentially the same. It is
based on the following expression for the number of ordered trees of order $N$ 
with
$n_k$ nodes having $k$ children:
$$
\frac{1}{N}\binom{N}{n_0,\ n_1,\ n_2 \ldots }
$$
if $n_1+2n_2+\ldots=N-1$, and $0$ otherwise
(see e.g.~Theorem 5.3.10 in \cite{Stanley:MR1676282}).

\bibliographystyle{alpha}
\bibliography{
LDP}
\end{document}